\documentclass[reqno,a4wide,english,11pt]{amsart}

\usepackage{amsmath,amssymb}
\usepackage{a4wide}

\usepackage[all]{xy}
\CompileMatrices
\newdir{ >}{{}*!/-10pt/\dir{>}}

\theoremstyle{plain}
\newtheorem{theorem}{Theorem}[section]

\theoremstyle{definition}

\newtheorem{emptythm}[theorem]{}

\newcommand{\too}{\longrightarrow}

\newcommand{\ul}[1]{\underline{#1}}

\newcommand{\BM}{\ul{\mathrm{h}}}

\DeclareMathOperator{\Hom}{Hom}
\DeclareMathOperator{\End}{End}

\DeclareMathOperator{\id}{id}

\DeclareMathOperator{\khar}{char}
\DeclareMathOperator{\Spec}{Spec}

\DeclareMathOperator{\res}{\mathrm{res}}

\newcommand{\Fields}{{\mathfrak Fields}}     

\newcommand{\SmF}[1]{{\mathbf Sm_{#1}}}

\newcommand{\Rings}{\mathfrak Rings}

\newcommand{\Chow}{\mathfrak{Chow}}
\newcommand{\Mot}{\mathfrak{Mot}}

\newcommand{\Tm}{\ul{\Z}}

\DeclareMathOperator{\CH}{\mathrm{CH}} 

\newcommand{\Cb}{\Omega} 

\newcommand{\Z}{\mathbb{Z}}

\newcommand{\Laz}{\mathbb{L}}

\newcommand{\cf}{\textsl{cf.}\ }
\newcommand{\eg}{\textsl{e.g.}\ }
\newcommand{\ie}{\textsl{i.e.}\ }

\newcommand{\locit}{\textsl{loc.cit.}\ }

\begin{document}

\title[Rost nilpotence and free theories]
{Rost nilpotence and free theories}

\author{Stefan Gille}
\email{gille@ualberta.ca}
\address{Department of Mathematical and Statistical Sciences,
University of Alberta, Edmonton T6G 2G1, Canada}

\author{Alexander Vishik}
\email{alexander.vishik@nottingham.ac.uk}
\address{School of Mathematical Sciences, University of Nottingham, Nottingham NG7 2RD,
United Kingdom}

\thanks{The work of the first author has been supported by an NSERC grant.}


\date{April 11, 2018}

\begin{abstract}
We introduce coherent cohomology theories~$\BM$ and prove that if such a theory
is moreover generically constant then the Rost nilpotence principle holds for projective
homogeneous varieties in the category of $\BM$-motives. Examples of such theories
are algebraic cobordism and its descendants the free theories.
\end{abstract}

\maketitle

\section{Introduction}
\label{IntroSect}\bigbreak

\noindent
Let~$X$ be a smooth projective quadric of dimension~$n$ over a field~$F$ of
characteristic~$\not= 2$. In~\cite{Ro98} Rost showed that a correspondence
$\alpha\in\CH_{n}(X\times_{F}X)$ which is rationally equivalent to zero over
a field extension~$E$ of~$F$ is nilpotent. A (purely) algebraic consequence
of this principle is that given a correspondence~$\pi$ of degree~$0$ on~$X$,
which becomes idempotent over an algebraic closure of~$F$, then~$\pi$ corresponds
to a direct summand of the quadratic. This implies in particular that there are no
{\it phantom} direct summands, \ie summands disappearing over an algebraic closure.
Therefore getting a complete decomposition of the (Chow-)motive of a projective quadric
is equivalent to finding a complete set of rational, \ie over the base field defined,
idempotent correspondences. Using this Rost gave in \locit a decomposition of the
Chow motive of a norm quadric which plays a crucial role in the proof of the Milnor
conjecture by Voevodsky~\cite{Voe03}.

\smallbreak

The usefulness of the nilpotence principle for quadrics in the category of Chow motives
is not limited to the proof of the Milnor conjecture. In fact, some of the recent advances in
the algebraic theory of quadratic forms rely on it, or its generalization to projective homogenous
varieties (as \eg orthogonal Grassmannians). A refinement of these 'geometric' methods for more
sophisticated cohomology theories than Chow theory should lead to further insights, but
requires the verification of the Rost nilpotence principle for quadrics and -- more general --
for projective homogeneous varieties for such theories. The aim of this work is to prove this
principle for free theories and projective homogeneous varieties. Moreover, we outline the general
context where the nilpotence principle holds. It is the context of the here introduced (generically)
constant {\it coherent} theories. Such a theory describes an environment where a cohomology functor
is defined not just for varieties of finite type over the ground field, but also for varieties over every
finitely generated extension of it, and moreover, these functors fit nicely with each other for different
extensions. A large supply of such theories are the so called {\it free theories}, \ie theories~$\BM_{\ast}$
which arise from algebraic cobordism by change of the group law: $\BM_{\ast}(X)=\BM_{\ast}(F)\otimes_{\Laz}\Cb_{\ast}(X)$,
where~$\Laz$ is the Lazard ring (cobordism of the point), and $\Cb_{\ast}$ algebraic cobrodism. But the class of {\it coherent}
theories is much bigger - we provide a construction which permits to produce such theories with very flexible properties. 
%
%
On our way we demonstrate also the failure of various aspects of the nilpotence principle outside the
{\it coherent} environment. This in particular shows that Rost nilpotence is not just a formal consequence 
of the axioms of oriented theories.

\bigbreak

\goodbreak
\section{Coherent theories}
\label{CohTheorySect}\bigbreak

\begin{emptythm}
\label{NotationConvSubSect}
{\bf Notations and conventions.}
Throughout this article fields are assumed to be of characteristic zero.

\smallbreak

For a field $F$ we denote by $\SmF{F}$ the category of smooth
quasi-projective schemes over $\Spec(F)$, and by $\Rings^{\ast}$ the category
of $\Z$-graded commutative rings.

\smallbreak

Given a field~$k$ we denote by $\Fields_{k}$ the category of finitely
generated field extensions of~$k$.

\smallbreak

If~$X$ is a scheme over a field~$F$ and~$L\supseteq F$ is a field
extension we set $X_{L}:=L\times_{F}X$. The function field of
an integral $F$-scheme~$X$ is denoted by~$F(X)$.
\end{emptythm}

\begin{emptythm}
\label{defCT}
{\bf Oriented cohomology- and Borel-Moore homology theories.}
We follow the definitions and conventions in the book~\cite{AC} on algebraic
cobordism by Levine and Morel, \ie an {\it oriented cohomology theory} over the field~$F$
is a functor
$$
\BM^{\ast}:\SmF{F}^{op}\rightarrow Rings^{\ast}
$$
with additional structure of push-forward maps $f_{\ast}:\BM^{\ast}(X)\too\BM^{\ast +d}(Y)$
along projective morphisms $f:X\too Y$ of (constant) relative codimension~$d$ subject to several
axioms including a projective bundle theorem, see~\cite[Def 1.1.2]{AC}.

\smallbreak

As usual we denote the pull-back $\BM^{\ast}(f)$ by~$f^{\ast}$ for
a morphism $f:X\too Y$ in~$\SmF{F}$.

\medbreak

If~$X\in\SmF{F}$ with connected components $X_{1},\ldots ,X_{l}$ one sets
$$
\BM_{\ast}(X)\, :=\;\bigoplus\limits_{i=1}^{l}\BM^{\dim X_{i} -\ast}(X_{i})\, .
$$
The assignment $X\mapsto\BM_{\ast}(X)$ defines then a Borel-Moore homology
theory on~$\SmF{F}$ in the sense of~\cite[Def.\ 5.1.3]{AC}, and by~\cite[Prop.\ 5.2.1]{AC}
$$
\BM_{\ast}\;\longleftrightarrow\,\BM^{\ast}
$$
is a one-to-one correspondence between oriented cohomology- and oriented Borel-Moore
homology theories on~$\SmF{F}$.

\smallbreak

Given an oriented Borel-Moore homology theory~$\BM_{\ast}$ on~$F$ we can extended
the theory to schemes of finite type over~$F$ by setting:
$$
\BM_{\ast}(Y)\, :=\;\operatorname{colim}_{V\rightarrow Y}\BM_{\ast}(V)\, ,
$$
where the limit runs over projective morphisms $V\rightarrow Y$ with~$V\in\SmF{F}$
with push-forward maps as transition maps.

\medbreak

In the following we consider theories with localization sequence (these correspond to
{\it oriented cohomology theories} of  \cite[Def. 2.1]{SU}). More precisely, we say an oriented
Borel-Moore homology theory~$\BM_{\ast}$, or equivalently the respective oriented cohomology
theory~$\BM^{\ast}$, over~$F$ has the {\it localization property} if given a smooth quasi-projective
$F$-scheme $X$ and a closed $F$-embedding $j:Z\rightarrow X$ with the open complement
$i:U\rightarrow X$ then there is an exact sequence
\begin{equation*}
\label{loc}
\tag{LOC}
\BM_{\ast}(Z)\stackrel{j_{\ast}}{\longrightarrow}\BM_{\ast}(X)\stackrel{i^{\ast}}{\longrightarrow}
\BM_{\ast}(U)\rightarrow 0\, .
\end{equation*}

\medbreak

\noindent
Then the same property will hold for any quasi-projective scheme $X$ over the base field~$F$,
see \cite[Sect.\ 2.2]{SU}.

\medbreak

\noindent
Examples of such theories are Chow groups, algebraic cobordism, and Grothendieck's~$K_{0}$
(in the latter example we ignore by some abuse of notation the grading).

\medbreak

An immediate consequence of the localization property
is the following useful lemma.
\end{emptythm}

\begin{emptythm}
\label{birOntoLem}
{\bf Lemma.}
{\it
Let~$\BM_{\ast}$ be an oriented Borel-Moore homology theory over the field~$F$
with localization property. Assume we have a projective morphism $\pi:X\too Y$ between
finite type $F$-schemes, whose restriction to $\pi^{-1}(U)$ is an isomorphism for some open
set~$U\subset Y$. Let $i:Z:=Y\setminus U\hookrightarrow Y$ be the closed
complement. Then
$$
\pi_{\ast}+i_{\ast}\, :\;\BM_{\ast}(X)\oplus\BM_{\ast}(Z)\,\too\,\BM_{\ast}(Y)
$$
is surjective.
}
\end{emptythm}

\begin{emptythm}
\label{defECT}
{\bf Definition.}
An {\it extended oriented cohomology theory} over the field~$k$ is a family of oriented
cohomology theories with localization property
$$
\BM_{F}^{\ast}\, :\;\SmF{F}\,\too\,\Rings^{\ast}\, ,\;\; F\in\Fields_{k}\, ,
$$
with the following additional data:

\smallbreak

Given $F\in\Fields_{k}$ and a smooth morphism $\rho:Y\too X$ in $\SmF{F}$ with $X$ integral
then there is a homomorphism of rings
$$
\BM_{F}^{\ast}(Y)\,\xrightarrow{\;\theta_{Y/X}\;}\,\BM_{F(X)}^{\ast}(F(X)\times_{X}Y)\, .
$$
(The map $\theta_{Y/X}$ might be interpreted as pull-back along
the upper row in the cartesian square
$$
\xymatrix{
F(X)\times_{X}Y \ar[r] \ar[d] & Y \ar[d]^-{\rho}
\\
\Spec F(X) \ar[r]_-{\iota} & X \rlap{\, ,}
}
$$
where $\iota:\Spec F(X)\too X$ is the generic point of~$X$.)

\bigbreak

These data are subject to the following four axioms:

\medbreak

\begin{itemize}
\item[(EC0)]
$\theta_{X/\Spec F}$ is the identity map for all $F\in\Fields_{k}$
and all $X\in\SmF{F}$.
\end{itemize}

\medbreak

Given a commutative diagram
$$
\xymatrix{
Y_{1} \ar[rr]^{f} \ar[rd] & & Y_{2} \ar[ld]
\\
 & X &
}
$$
in~$\SmF{F}$ with~$X$ integral for some finitely generated field extension~$F$ of~$k$ we denote by
$\tilde{f}:F(X)\times_{X}Y_{1}\too F(X)\times_{X}Y_{2}$ the induced morphism on the generic
fibers. Then:

\medbreak

\begin{itemize}
\item[(EC1)]
For any morpism~$f$ we have
$$
\theta_{Y_{1}/X}\circ f^{\ast}\, =\,\tilde{f}^{\ast}\circ\theta_{Y_{2}/X}\, ;
$$

\smallbreak

\item[(EC2)]
If~$f$ is projective of constant relative dimension
$$
\tilde{f}_{\ast}\circ\theta_{Y_{1}/X}\, =\,\theta_{Y_{2}/X}\circ f_{\ast}\, .
$$
\end{itemize}

\medbreak

Given a diagram of cartesian squares
$$
\xymatrix{
{\widetilde{Z}} \ar[r] \ar[d] & F(X)\times_{X}Z \ar[r] \ar[d]_-{{\tilde{\sigma}}} & Z \ar[d]^-{\sigma}
\\
{\Spec F(Y)} \ar[r] & F(X)\times_{X}Y \ar[d] \ar[r] & Y \ar[d]^-{\rho}
\\
 & {\Spec F(X)} \ar[r] & X
}
$$
with~$X,Y$ integral and $\rho ,\sigma$ morphisms in~$\SmF{F}$ for some $F\in\Fields_{k}$
(note that $F(Y)$ and the function field of $F(X)\times_{X}Y$ coincide). Then:

\medbreak

\begin{itemize}
\item[(EC3)]
In this situation the following holds:
$$
\theta_{(F(X)\times_{X}Z)/(F(X)\times_{X}Y)}\circ\theta_{Z/X}=\theta_{Z/Y}\, :\; A_{F}^{\ast}(Z)\;\too\;
A_{F(Y)}^{\ast}(\widetilde{Z})\, .
$$
\end{itemize}
\end{emptythm}

\begin{emptythm}
\label{relCohGr}
{\bf The group $\BM^{\ast}_{F}(Y/X)$.}
Let~$\BM^{\ast}$ be an extended oriented cohomology theory over the field~$k$
and~$F$ a finitely generated field extension of~$k$.

\smallbreak

For a smooth morphism $Y\rightarrow X$ in $\SmF{F}$ with $X$ integral we set
$$
\BM^{\ast}_{F}(Y/X)\, :=\;\operatorname{colim}_{U\rightarrow X}\BM^{\ast}_{F}(U\times_{X}Y)\, ,
$$
where the colimit runs over all open subschemes~$U$ of $X$.

\smallbreak

Let $U\subseteq X$ be an open subscheme. Then we have a homomorphism
$$
\theta_{U\times_{X}Y/X}\, :\;\BM^{\ast}_{F}(U\times_{X}Y)\,\too\,\BM^{\ast}_{F(X)}(\widetilde{Y})\, ,
$$
where we used that $\widetilde{Y}:=F(X)\times_{X}Y=F(X)\times_{X}(U\times_{X}Y)$.
By axiom~(EC1) we have $\theta_{Y/X}=\theta_{U\times_{X}Y/X}\circ q^{\ast}$, where
$q:U\times_{X}Y\too Y$ is the projection, and more generally
$$
\theta_{U\times_{X}Y/X}\, =\,\theta_{V\times_{X}Y/X}\circ (i\times\id_{Y})^{\ast}
$$
if $i:V\hookrightarrow U$ is another open subscheme of~$X$ contained in~$U$.
Hence we have a canonical map
$$
\varphi_{Y/X}\, :\; \BM^{\ast}_{F}(Y/X)\,\too\,\BM^{\ast}_{F(X)}(\widetilde{Y})\, .
$$
\end{emptythm}

\begin{emptythm}
\label{DefCohCT}
{\bf Definition.}
An extended oriented cohomology theory over $k$ is called {\it coherent} if the map
$\varphi_{Y/X}$ is an isomorphism for all smooth $F$-morphsims $Y\too X$ and all $F\in\Fields_{k}$.
\end{emptythm}

\begin{emptythm}
\label{ECT-standardExpl}
{\bf Examples.}
\begin{itemize}
\item[(i)]
Let $\Cb^{\ast}$ be algebraic cobordism as defined by Levine and Morel~\cite{AC},
$\rho:Y\too X$ a morphism in~$\SmF{F}$ with~$X$ integral for some $F\in\Fields_{k}$.
Set $E:=F(X)$ and~$\widetilde{Y}:=E\times_{X}Y$. Given a cobordism cycle
$$
[g:W\too Y\, ,\; L_{1},\ldots ,L_{r}]
$$
over~$Y$ then the pull-back
$$
\big[\, \Spec E\times_{X}W\too\widetilde{Y}\, ,\; p^{\ast}_{W}(L_{1}),\ldots ,
p^{\ast}_{W}(L_{r})\,\big]\, ,
$$
where $p_{W}:\Spec E\times_{X}W\too W$ is the projection, is a cobordism cycle
over~$\widetilde{Y}$. It is straightforward to check that this induces a homomorphism
$\theta_{Y/X}:\Cb^{\ast}(Y)\too\Cb^{\ast}(\widetilde{Y})$. We will show below, see
Corollary~\ref{Omega-EOCT}, that this map obeys the axioms (EC1)-(EC3), and so
algebraic cobordism is an extended oriented cohomology theory over all fields of
characteristic~$0$.

\smallbreak

\item[(ii)]
From the fact that algebraic cobordism is an extended oriented cohomology theory it
follows also that Chow groups with rational or algebraic equivalence and arbitrary coefficients,
and also Grothendieck's~$K_{0}$ are extended oriented cohomology theories. However this can
be seen more easily directly.
\end{itemize}
\end{emptythm}

\begin{emptythm}
\label{Exa-notCoh}
{\bf Remark.}
Not all extended oriented cohomology theories are {\it coherent}. An example is provided by
$\CH^{\ast}_{alg}$. Indeed, consider the smooth projection $C\times C\rightarrow C$, where $C$
is an elliptic curve over $k$. Then in $\CH^{\ast}_{alg}(C\times C)$ the class of the diagonal $[\Delta]$
is not equal to the linear combination of $[C\times p]$ and $[p\times C]$, for a $k$-rational point $p$,
since it is not so even modulo homological equivalence. Hence, from the localization sequence (and the
fact that classes of all rational points on $C$ are equivalent modulo algebraic equivalence) we obtain
that the restriction of $[\Delta]$ to $\left(\CH^{\ast}_{alg}\right)_{k}(C\times C/C)$ is different from that
of $[C\times p]$. On the other hand, the restriction of these two classes to $\left(\CH^{\ast}_{alg}\right)_{k(C)}(C_{k(C)})$
are just classes of two rational points which are equal. Hence, the map $\varphi_{C\times C/C}$ is not injective.

\medbreak

But below we will show that all the theories obtained from $\Cb^{\ast}$ by change of coefficients
(the so-called {\it free theories}) are coherent.
\end{emptythm}

\begin{emptythm}
\label{MainLemma}
{\bf Lemma.}
{\it
Let $\BM^{\ast}:=\Cb^{\ast}$ be algebraic cobordism and $F$~a field. Assume we have a
commutative diagram
$$
\xymatrix{
Y \ar[rr]^-{f} \ar[rd]_-{p_{Y}} & & X \ar[ld]^-{p_{X}}
\\
 & U &
}
$$
in~$\SmF{F}$ with~$U$ integral and $p_{X},p_{Y}$ smooth. Denote by~$\widetilde{X}$ and~$\widetilde{Y}$
the generic fiber of~$p_{X}$ and~$p_{Y}$, respectively, and let $\tilde{f}:\widetilde{Y}\too\widetilde{X}$ be
the morphism induced by~$f$. Then we have:

\medbreak

\begin{itemize}
\item[(i)]
$\theta_{X/U}(\alpha\cdot\beta)=\theta_{X/U}(\alpha)\cdot\theta_{X/U}(\beta)$
for all $\alpha,\beta\in\Cb^{\ast}(X)$;

\bigbreak

\item[(ii)]
$\widetilde{f}^{\ast}(\theta_{X/U}(\alpha))=\theta_{Y/U}(f^{\ast}(\alpha))$ for all $\alpha\in\Cb^{\ast}(X)$;
and

\bigbreak

\item[(iii)]
if~$f$ is projective, then
$\widetilde{f}_{\ast}(\theta_{Y/U}(\gamma))=\theta_{X/U}(f_{\ast}(\gamma))$ for all
$\gamma\in\Cb^{\ast}(Y)$.
\end{itemize}
}

\smallbreak

\begin{proof}
The proof uses the following property of algebraic cobordism. Assume we have
a closed regular embedding $j: V\hookrightarrow W$ in~$\SmF{F}$ for some field~$F$.
Then it is shown in Levine and Morel~\cite[Prop.\ 3.3.1]{AC} that~$\Cb^{\ast}(W)$ is
generated by cobordism cycles $(h:Z\too W)$ with~$h$ transversal to~$j$. The image
of such a cycle under~$j^{\ast}$ is then just the pull-back of cycles. In fact, since~$V,W$
are smooth over~$F$ the fiber product $Z\times_{W}V$ is regular and so smooth, as $\khar F=0$.
Hence we can apply~\cite[Cor.\ 6.5.5]{AC} and get
$$
j^{\ast}\big(\, [h:Z\too W]\,\big)\, =\, [h':Z\times_{W}V\too V]\, ,
$$
where~$h'$ is the projection.

\smallbreak

This implies~(ii) if~$f$ is a regular embedding. For the general case of~(ii) observe that~$f$
factors $Y\xrightarrow{f\times\id}X\times_{F}Y\rightarrow X$ and so, since~$X$ and~$Y$ are smooth,
into a composition of a regular embedding and a smooth morphism. This reduces then the proof of~(ii)
to the case of a smooth equi-dimensional morphism of constant relative dimension, where it follows
from the definition of the pull-back in this case, see~\cite[Sect.\ 2.1]{AC}. At the same time, the assertion~(iii)
is an immediate consequence of the definition of the push-forward of cobordism cycles along
projective maps.

\medbreak

Finally we prove~(i). Let $\Delta_{X/U}:X\too X\times_{U}X$ and
$\Delta_{X/F}:X\too X\times_{F}X$ be the respective diagonal maps.
Consider $\Delta_{X/U}$ as a $U$-morphisms and denote $\widetilde{\Delta}_{X/U}$
the induced map between the generic fibers. Now $F(U)\times_{U}(X\times_{U}X)$ is
naturally isomorphic to the fiber product of $\widetilde{X}=F(U)\times_{U}X$ with itself
over~$F(U)$, and under this identification~$\widetilde{\Delta}_{X/U}$ corresponds to the
diagonal map $\Delta_{\widetilde{X}/F(U)}:\widetilde{X}\too\widetilde{X}\times_{F(U)}\widetilde{X}$.
Hence by~(ii) we have
\begin{equation}
\label{MainLemPfEq1}
\theta_{X/U}\circ\Delta_{X/U}^{\ast}=\Delta_{\widetilde{X}/F(U)}^{\ast}\circ\theta_{X\times_U X/U}\, .
\end{equation}
By definition of the product in algebraic cobordism we have
$$
\alpha\cdot\beta\, =\,\Delta^{\ast}_{X/F}(\alpha\times\beta)\quad\mbox{and}\quad
\theta_{X/U}(\alpha)\cdot\theta_{X/U}(\beta)\, =\,\Delta^{\ast}_{\widetilde{X}/F(U)}
\big(\theta_{X/U}(\alpha)\times\theta_{X/U}(\beta)\big)\, .
$$
Hence by~(\ref{MainLemPfEq1}) it remains to show
\begin{equation}
\label{MainLemPfEq2}
\theta_{X\times_{U} X/U}\big(\delta_{X/U}^{\ast}(\alpha\times\beta)\big)\, =\,
\theta_{X/U}(\alpha)\times\theta_{X/U}(\beta)
\end{equation}
for all $\alpha,\beta\in\Cb^{\ast}(X)$, where $\delta_{X/U}$ is the regular immersion
$X\times_{U}X\too X\times_{F}X$.

\medbreak

If $\rho:U'\rightarrow U$ is an open embedding and $\chi\, :\; X':=U'\times_{U}X\too X$
the projection then we have
$$
\delta_{X'/U'}^{\ast}\circ (\chi\times\chi)^{\ast}\, =\, (\chi\times_{\rho}\chi)^{\ast}\circ\delta_{X/U}^{\ast}\quad
\mbox{and}\quad \theta_{X\times_{U} X/U}=\theta_{X'\times_{U'}X'/U'}\circ (\chi\times_{\rho}\chi)^{\ast}.
$$
Thus to prove~(\ref{MainLemPfEq2}) we can replace $U$ by an open subscheme of~$U$.

\smallbreak

After these reductions and remarks we can finish the proof of~(i). We can assume that~$\alpha$
is represented by $(v:V\too X)$ and~$\beta$ by $(w:W\too X)$ for some $V,W\in\SmF{F}$. The
generic fiber $F(U)\times_{U}W$ is (as localization of a regular scheme) also regular and therefore
smooth over~$F(U)$ since (by our general assumption) $\khar F=0$. Therefore replacing~$U$ by an
open subscheme of it we can assume that the composition of morphisms $W\too X\xrightarrow{p_{X}}U$
is smooth. Then $V\times_{F}W\too X\times_{F}X$ is transversal to the regular embedding
$\delta_{X/U}:X\times_{U}X\too X\times_{F}X$ and $V\times_{U}W$ is smooth over~$F$. Hence
by~\cite[Cor.\ 6.5.5]{AC} the pull-back  $\delta^{\ast}_{X/U}(\alpha\times\beta)$ is represented by the
cartesian product $(V\times_{U} W\rightarrow X\times_{U} X)$, and so we have
$\theta_{X\times_U X/U}=\theta_{X/U}(\alpha)\times\theta_{X/U}(\beta)$.
The claim follows.
\end{proof}
\end{emptythm}

\begin{emptythm}
\label{Omega-EOCT}
{\bf Corollary.}
{\it
Algebraic cobordism (of Levine and Morel) is an extended oriented cohomology theory.
}

\begin{proof}
The axiom~(EC0) is obvious, and by Lemma \ref{MainLemma} above
$\theta_{Y/X}$ is a homomorphism of rings satisfying also~(EC1) and~(EC2).

\smallbreak

For~(EC3): It is enough to check this on 'standard' cobordism cycles
$(w: W\too X)$, where it is straightforward and amounts to the identification
of two different ways to express the generic fiber of a map to a fibration.
\end{proof}

\smallbreak

For the proof of the following theorem we have to use the generators and relations of
algebraic cobordism.
\end{emptythm}

\begin{emptythm}
\label{Omega-Coh}
{\bf Theorem.}
{\it
Algebraic cobordism is a coherent theory.
}

\begin{proof}
Let $F$ be a field of characteristic zero, and $Y\stackrel{\pi}{\longrightarrow} X$ be a
smooth morphism in $\SmF{F}$ with $X$ integral, $E=F(X)$ and $\widetilde{Y}$ the
generic fiber. Let $\Cb^{*}_{F}(Y/X)=\operatorname{colim}_{U\rightarrow X}\Cb^{*}(Y\times_{X}U)$
over all open subschemes of $X$ and
$$
\varphi_{Y/X}:\Cb^{\ast}_{F}(Y/X)\too\Cb^{\ast}_{E}(\widetilde{Y})
$$
be the canonical map. We need to show that $\varphi_{Y/X}$ is an isomorphism.

Combining the definitions of $\Omega^{\ast}$ from \cite{AC} and \cite{LePa09}
we have an exact sequence:
$$
{\mathcal R}_{\ast}(Q)\rightarrow {\mathcal Z}_{\ast}(Q)\rightarrow \Cb_{*}(Q)\rightarrow 0.
$$
Here ${\mathcal Z}_{\ast}(Q)$ is a free abelian group generated by isomorphism classes of
{\it cobordism cycles} $[v:V\rightarrow Q,L_1,\ldots,L_r]$, where $v$ is projective, $V$ is smooth, and
$L_{1},\ldots,L_{r}$ are (non-ordered) line bundles on $V$ -
see \cite[Def 2.1.6]{AC}.

The relations ${\mathcal R}_{\ast}(Q)$ form a free abelian group which splits as a direct sum
$$
{\mathcal R}_{\ast}(Q)\, =\,{\mathcal R}_{\ast}^{Sect}\oplus {\mathcal R}_{\ast}^{Dim}\oplus {\mathcal R}_{\ast}^{MPR}\, .
$$
The three groups on the right hand side and their maps to~${\mathcal Z}_{\ast}(Q)$ are given by:
\begin{itemize}
\item
${\mathcal R}_{\ast}^{Dim}(Q)$ is the free abelian group generated by isomorphism classes of cobordism cycles
$[v:V\rightarrow Q,L_1,\ldots,L_r]$ with $r>\dim V$ with the obvious map to ${\mathcal Z}_{\ast}(Q)$.

\medbreak

\item
${\mathcal R}_{\ast}^{Sect}(Q)$ is the free abelian group generated by the isomorphism classes of pairs
$[Z\stackrel{i}{\subset} V\stackrel{v}{\rightarrow}Q,L_{1},\ldots,L_{r}]$ of a cobordism cycle plus a smooth divisor
on~$V$. Such a pair is mapped to the difference
$$
[V\stackrel{v}{\rightarrow}Q,L_{1},\ldots,L_{r},O(Z)]-
[Z\stackrel{v\circ i}{\rightarrow}Q,i^{\ast}L_{1},\ldots,i^{\ast}L_{r}]
$$
in ${\mathcal Z}_{\ast}(Q)$.

\medbreak

\item
Finally, ${\mathcal R}_{\ast}^{MPR}$ is the free abelian group generated by the isomorphism classes of
{\it multiple point relations} $[w:W\rightarrow Q\times{\mathbb{P}}^{1},M_{1},\ldots,M_{s}]$, where $w$ is projective,
$W$ is smooth, $W_{0}=w^{-1}(Q\times 0)$ and $W_{1}=w^{-1}(Q\times 1)$ are divisors with strict normal crossings
on~$W$, see~\cite[Def 3.1.4]{AC}), and $M_{1},\ldots, M_{s}$ are line bundles on~$W$.

\smallbreak

Such an element is mapped to $w_{\ast}([W_{0}\rightarrow W,M_{1},\ldots, M_{s}]-[W_{1}\rightarrow W,M_{1},\ldots, M_{s}])$,
where $[W_l\rightarrow W]$ is the {\it combinatorial divisor class}, see \cite[Def 3.15]{AC}, and where we lift the coefficients of
the universal formal group law to cobordism cycles from ${\mathcal Z}_{\ast}(\Spec(k))$. Recall, that if $D=\sum_{i=1}^{m}n_{i}\cdot D_{i}$
is a divisor with strict normal crossing on $W$, with $L_{i}=O(D_i)$ and $d_{J}:D_{J}=\cap_{j\in J}D_{j}\rightarrow W$ for
$J\subset\{1,\ldots,m\}$, the faces of our divisor, we can write down the formal sum
$$
[n_{1}]\cdot_{F}u_{1}+_{F}\ldots+_{F}[n_m]\cdot_{F}u_{m}
\hspace{5mm}\text{as}\hspace{2mm}
\sum_{J\subset\{1,\ldots ,m\}}u^{J}\cdot F_{J}^{n_{1},\ldots ,n_{m}}(u_{1},\ldots ,u_{m})\, ,
$$
where
$$
F_{J}^{n_{1},\ldots ,n_{m}}(u_{1},\ldots ,u_{m})\;\in\,\Laz[[u_{1},\ldots ,u_{m}]]
$$
are some power series with $\Laz$-coefficients, which can be chosen canonically.
Then the combinatorial divisor class $[D\rightarrow W]$ is just the sum
$\sum_{J\subset\{1,\ldots,m\}}[d_{J}:D_{J}\rightarrow W,F_{J}^{n_{1},\ldots ,n_{m}}(L^{J}_{1},\ldots, L^{J}_{m})]$, where
$L^{J}_{i}=d_J^{*}L_{i}$. Here we ignore all the terms of degree $>\dim W$ (these are covered by ${\mathcal R}_{\ast}^{Dim}$), and so,
it is a finite sum. Note, that ${\mathcal R}_{\ast}^{MPR}$ contains the {\it geometric cobordism relations}, \cf~\cite[Def 2.3.2]{AC},
as well as {\it double point relations}, see \cite{LePa09}, and so, the (FGL) relations of \cite{AC}.
\end{itemize}

\smallbreak

\noindent
Finally, the map ${\mathcal Z}_{\ast}(Q)\rightarrow\Cb_{\ast}(Q)$ is given by
$$
[v:V\rightarrow Q,L_{1},\ldots,L_{r}]\mapsto v_{*}(c_{1}(L_{1})\cdot\ldots\cdot c_1(L_{r})).
$$
It follows from \cite[Thm 2.3.13]{AC} that in ${\mathcal R}_{\ast}^{Dim}$ we can take any number
$N\geq\dim V$ instead of $\dim V$ (modifying accordingly the interpretation of ${\mathcal R}_{\ast}^{MPR}$).

\medbreak

We have a commutative diagram with exact rows:
\begin{equation*}
\def\objectstyle{\displaystyle}
\xymatrix{
\tilde{{\mathcal R}}_{\ast}(\widetilde{Y}) \ar[r] & {\mathcal Z}_{\ast}(\widetilde{Y}) \ar[r] & \Cb_{*}^{E}(\widetilde{Y}) \ar[r] & 0\\
\operatornamewithlimits{colim}_{U\rightarrow X}{\mathcal R}_{\ast}(Y\times_X U) \ar[r] \ar[u]_-{\rho_{Y/X}} &
\operatornamewithlimits{colim}_{U\rightarrow X}{\mathcal Z}_{\ast}(Y\times_X U) \ar[r] \ar[u]_-{\psi_{Y/X}} &
\Cb_{*}^{F}(Y/X) \ar[r] \ar[u]_-{\varphi_{Y/X}} &  0,
}
\end{equation*}
where in the definition of $\tilde{{\mathcal{R}}}_{\ast}^{Dim}(\widetilde{Y})$ we swap
$\dim \widetilde{V}$ by $\dim \widetilde{V}+\dim X$.

\smallbreak

We show now that $\rho_{Y/X}$ is onto and $\psi_{Y/X}$ is an isomorphism, from which it
follows that~$\varphi_{Y/X}$ is an isomorphisms, too, hence finishing the proof.

\medbreak

We start showing~$\psi_{Y/X}$ is onto. In fact, if we have a projective map
$\widetilde{v}:\widetilde{V}\rightarrow\widetilde{Y}$ in $\SmF{E}$, we can always extend
it to some projective map $v:V\rightarrow Y$. Moreover, using the resolution of singularities,
we can make $V$ smooth. Also, any line bundle can be extended from $\widetilde{V}$ to $V$.
This shows that $\psi_{Y/X}$ is surjective.

\smallbreak

If two maps $v:V\rightarrow Y\times_{X} U$ and $v':V'\rightarrow Y\times_{X} U'$ are isomorphic
over the generic fiber $\widetilde{Y}$, then this isomorphism is actually defined over some open
neighborhood $U''\subset U\cap U'$ of the $\Spec(E)$. And if two line bundles $L$ and $L'$ on $V$ are
isomorphic when restricted to $\widetilde{V}$, then these differ in the Picard group by a class
of a line bundle coming from $U$ (since $\widetilde{V}$ is obtained from $V$ by removing preimages
of some divisors from $U$). So, these become isomorphic over the preimage of some neighborhood of $\Spec(E)$.
Hence, the map $\psi_{Y/X}$ is injective, and so an isomorphism.

\smallbreak

By the same arguments, the modified relations $\tilde{{\mathcal R}}_{\ast}^{Dim}(\widetilde{Y})$ are covered by
${\mathcal R}_{\ast}^{Dim}(Y)$.

\smallbreak

If we have some projective map $\widetilde{v}:\widetilde{V}\rightarrow\widetilde{Y}$ with the smooth divisor
$\widetilde{Z}\subset\widetilde{V}$, as above, using the results of Hironaka \cite{Hi}, we can extend $\widetilde{v}$ to a
projective map $v:V\rightarrow Y$ from a smooth $V$, and can extend $\widetilde{Z}$ to a smooth divisor $Z$ on $V$.
Also, we can extend any line bundle from $\widetilde{V}$ to $V$. Hence, the map $\rho_{Y/X}$ is surjective on
${\mathcal R}_{\ast}^{Sect}$.

\smallbreak

If we have some multiple point relation $\widetilde{w}:\widetilde{W}\rightarrow\widetilde{Y}\times{\Bbb{P}}^1$,
then we can extend $\widetilde{w}$ to a projective map $w:W\rightarrow Y\times{\mathbb{P}}^1$. Using the results of
Hironaka~\cite{Hi} we can resolve the singularities of $W$ and make the pre-images of $Y\times 0$ and
$Y\times 1$ divisors with strict normal crossing without changing the generic fiber $\widetilde{w}$.
Thus, the relations ${\mathcal R}_{\ast}^{MPR}(\widetilde{Y})$ are covered by $\rho_{Y/X}$.
\end{proof}
\end{emptythm}

\begin{emptythm}
\label{free-theories}
{\bf Free theories.}
Any formal group law ${\mathcal F}$ over a ring~$R$ (possibly graded)
is induced from the universal one $(\Laz_{\ast},{\mathcal F}_U)$ via the unique homomorphism of rings
$\Laz_{\ast}\too R$ and one can assign to it a theory
$$
X\,\longmapsto\,\BM_{\mathcal F}(X)\, :=\; R\otimes_{\Laz}\Cb_{\ast}(X)\, .
$$
Such a theory is called {\it free} (see \cite[Rem 2.4.14]{AC}) and inherits from algebraic cobordism
all basic properties of an oriented cohomology theory.
Examples of free theories are cobordism itself with arbitrary coefficients, Chow groups
with rational equivalence and arbitrary coefficients, $K_0$, and Morava $K$-theories.

\smallbreak

Since all our structural maps in $\Cb_{\ast}$ are $\Laz$-linear, it follows that a free
theory~$\BM_{\mathcal F}$ is an {\it extended oriented cohomology theory}.
And since colimits commutes with tensor products we obtain:

\medbreak

\noindent
{\bf Corollary.}
{\it Free theories are coherent.}
\end{emptythm}

\medbreak

\begin{emptythm}
\label{BaseChangeSubSect}
{\bf Base change.}
Let~$\BM^{\ast}$ be an extended oriented cohomology theory over the field~$k$.

\smallbreak

Given~$L\supseteq F$ in~$\Fields_{k}$ there exists a smooth integral $F$-scheme~$U$,
such that $L=F(U)$. Let $\iota:\Spec L=\Spec F(U)\hookrightarrow U$ be the embedding of the generic point.

\smallbreak

Let now~$X\in\SmF{F}$. We have a smooth morphism $U\times_{F} X\rightarrow U$, and we will denote the
corresponding map $\theta_{U\times X/U}$ by
$$
(\iota\times\id_{X})^{\ast}\, :\;\BM_{i}(U\times_{F}X)\,\too\,\BM_{i-\dim U}(X_{F(U)})\,
=\,\BM_{i-\dim U}(X_{L})\, .
$$
Composing it with the map $\pi^{\ast}$, where $\pi:U\times_{F}X\rightarrow X$ is the projection we obtain
the {\it base change homomorphism}:
$$
\BM_{\ast}(X)\,\too\,\BM_{\ast}(X_{L})\, ,\;\alpha\,\longmapsto\,\alpha_{L}\, ,
$$
which we denote by~$\res_{L/F}$. This homomorphism does not depend on the choice
the model~$U$. In fact, if~$U'$ is another model then~$U'$ is birational to~$U$, and so
by~(EC1) we are reduced to show that if~$j:U'\hookrightarrow U$ is an open subscheme
then $\theta_{U\times X/U}=\theta_{U'\times X/U'}\circ (j\times\id_{X})^{\ast}$. This can be
seen as follows: Applying~(EC3) to the sequence of morphisms
$$
U'\times_{F}X\,\xrightarrow{\; p\;}\, U'\,\xrightarrow{\: j\:}\, U\, ,
$$
where~$p$ is the projection, and taking~(EC0) into account we obtain
$\theta_{U'\times X/U'}=\theta_{U'\times X/U}$. From this equality we get
$\theta_{U\times X/U}=\theta_{U'\times X/U'}\circ (j\times\id_{X})^{\ast}$ by~(EC1).

\smallbreak

From the fact that $\theta_{X\times U/U}$ and $\pi^*$ are ring homomorphisms and
from~(EC1) and~(EC2) we get:
\end{emptythm}

\begin{emptythm}
\label{NatPropLem}
{\bf Lemma.}
{\it
Let~$\BM_{\ast}$ a theory over a field~$k$ be as above, $F\in\Fields_{k}$, $X$ and~$Y$ smooth $F$-schemes, and
$L=F(U)\supseteq F$ a field extension, where~$U$ is a smooth and integral $F$-scheme.
Denote~$\iota$ the inclusion of the generic point of~$U$ into~$U$. Let further $f,g:Y\too X$ be a morphism,
respectively, projective morphism in~$\SmF{F}$.

\medbreak

Then:

\smallbreak

\begin{itemize}
\item[(i)]
$(\iota\times\id_{X})^{\ast}(\alpha\cdot\beta)=
(\iota\times\id_{X})^{\ast}(\alpha)\,\cdot\, (\iota\times\id_{X})^{\ast}(\beta)$
for all $\alpha,\beta\in\BM_{\ast}(U\times X)$.

\bigbreak

\item[(ii)]
$f_{L}^{\ast}\circ (\iota\times\id_{X})^{\ast}=(\iota\times\id_{Y})^{\ast}\circ (\id_{U}\times f)^{\ast}$.

\bigbreak

\item[(iii)]
$g_{L\,\ast}\circ (\iota\times\id_{Y})^{\ast}=(\iota\times\id_{X})^{\ast}\circ (\id_{U}\times g)_{\ast}$.

\bigbreak

\item[(iv)]
$(\alpha\cdot\beta)_{L}=\alpha_{L}\cdot\beta_{L}$
for all $\alpha,\beta\in\BM_{\ast}(X)$.
\end{itemize}
}
\end{emptythm}

\goodbreak
\section{$\BM$-motives}
\label{hMotivesSect}\bigbreak

\begin{emptythm}
\label{deg0CorrSubSect}
{\bf $\BM$-correspondences of degree~$0$.}
The following definitions and assertions are an adaption of Manin's~\cite{Ma68} for Chow groups
(see also Fulton~\cite[Chap.\ 16]{INT}). We recall them here for the sake of completeness.

\medbreak

Let~$\BM_{\ast}$ be an oriented Borel-Moore homology theory over a field~$k$ and~$X,Y$
smooth projective $k$-schemes. Denote $X_{1},\ldots ,X_{l}$ the connected
components of~$X$. An {\it $\BM$-correspondence} between~$X$ and~$Y$ is an
element of $\BM_{\ast}(X\times_{k}Y)$, and an {\it $\BM$-correspondence of degree~$0$}
between~$X$ and~$Y$ is an element
$$
\alpha\,\in\;\bigoplus\limits_{i=1}^{l}\BM_{\dim X_{i}}(X_{i}\times_{k}Y)\;\subseteq
\BM_{\ast}(X\times_{k}Y)\, .
$$
We indicate the latter situation by $\alpha:X\rightsquigarrow Y$.

\smallbreak

The {\it composition} of two correspondences $\alpha: X\rightsquigarrow Y$ and
$\beta:Y\rightsquigarrow Z$ is defined for connected smooth and projective
schemes~$X$ and~$Y$ as follows:
$$
\beta\circ\alpha\, :=\; p_{XZ\,\ast}\big(p^{\ast}_{XY}(\alpha)\cdot p^{\ast}_{YZ}(\beta)\big)\, :\;
X\,\rightsquigarrow\, Z\, ,
$$
where~$\cdot$ denotes the product in $\BM_{\ast}(X\times_{k}Y\times_{k}Z)$ and
$p_{XY}$, $p_{YZ}$, and $p_{XZ}$ denote the indicated projections from
$X\times_{k}Y\times_{k}Z$ to $X\times_{k}Y$, $Y\times_{k}Z$, and $X\times_{k}Z$,
respectively. This composition extends then linearly to non connected schemes,
and is associative with the image of the $1\in\BM_{\ast}(X)$ under the push-forward
along the diagonal embedding $X\too X\times_{k}X$ acting as identity.

\smallbreak

We denote the set of $\BM$-correspondences $X\rightsquigarrow Y$ of degree~$0$
by $\Hom_{k}(X,Y)_{\BM}$, respectively by $\End_{k}(X)_{\BM}$ if~$X=Y$.

\medbreak

Given an $\BM$-correspondence $\alpha:X\rightsquigarrow X$ of
degree~$0$, where~$X$ is a smooth and projective $k$-scheme,
it acts on $\BM_{i}(X)$ by
$$
\alpha_{\ast}(\gamma)\, :=\;\alpha\circ\gamma\, .
$$
\end{emptythm}

\begin{emptythm}
\label{MotSubSect}
{\bf $\BM$-motives.}
It follows from Manin's~\cite{Ma68} that the category whose objects are the smooth
projective $k$-schemes and whose groups of morphisms are given by $\Hom_{k}(X,Y)_{\BM}$
is an additive category, called the category of $\BM$-correspondences of degree~$0$. The
idempotent completion of this category is the {\it category of $\BM$-motives over~$k$}, which
we will denote $\Mot_{\BM}(k)$. If $\BM_{\ast}=\CH_{\ast}$ this is the well known category of
effective Chow motives~$\Chow (k)$.

\smallbreak

In this category the projective line decomposes into two summands, one isomorphic to the
motive of the point~$\Spec k$, which is denoted by $\Tm_{\BM}$, and the other summand
is called the {\it Tate-} or {\it Lefschetz motive} and denoted by $\Tm_{\BM}\{1\}$.

\smallbreak

The cartesian product induces a 'tensor' product on the motives. In particular
we have the (positive) {\it Tate twists}
$$
\Tm_{\BM}\{r\}\, :=\;\Tm_{\BM}\{1\}^{\otimes\, r}\qquad\mbox{and}\qquad
X\{r\}\, :=\,\Tm_{\BM}\{r\}\otimes X
$$
for~$X$ a smooth and projective $k$-scheme and~$r\geq 0$ an integer.
\end{emptythm}

\begin{emptythm}
\label{RNprincipleSubSect}
{\bf The Rost nilpotence principle.}
Let now~$\BM_{\ast}$ be an extended oriented cohomology theory.
By some abuse of notation we suppress the lower supscript in~$\BM_{\ast\, E}$
and write $\BM_{\ast}$ instead.

\medbreak

We say that a smooth projective variety~$X$ over the field~$k$
satisfies the {\it Rost nilpotence principle} for the theory~$\BM_{\ast}$ if for all
field extensions $F\supseteq k$ the kernel of the base change map
$$
\res_{F/k}\, :\;\End_{k}(X)_{\BM}\,\too\,\End_{F}(F\times_{k}X)_{\BM}\, ,\;\alpha\,\longmapsto\,\alpha_{F}
$$
consists of nilpotent correspondences, \ie correspondences~$\alpha$, such that
$\alpha^{\circ\, N}=0$ for some natural number~$N$.

\bigbreak

For~$\BM_{\ast}=\CH_{\ast}$ it is conjectured that the Rost nilpotence principle
holds for all smooth projective schemes for all fields. This has been proven in some cases,
for instance quadrics by Rost~\cite{Ro98} (from where the name of the principle comes from) and independently by the
second named author~\cite{IMQ},
projective homogeneous varieties by Chernousov, the first named author, and Merkurjev~\cite{ChGiMe05}
and independently by Brosnan~\cite{Bro05}; and for surfaces over fields of characteristic~$0$
in~\cite{Gi10,Gi14}.

\smallbreak

The arguments in the article~\cite{ViYa07} of Yagita and the second named author imply that if~$X$
is a smooth projective scheme over a field of characteristic zero and Rost nilpotence is known for~$X$
in the category of Chow motives then it is known for~$X$ in the category of $\Cb$-motives. In particular
it is shown there that a decomposition in the category of Chow motives implies a decomposition in the
category of $\Cb$-motives and so also in the category of $\BM$-motives for all oriented cohomology theories~$\BM_{\ast}$
by the universality of algebraic cobordism, see~\cite[Chap.\ 7]{AC}.

\smallbreak

As by the main result of~\cite{ChGiMe05}, see also Brosnan~\cite{Bro05}, an isotropic projective homogeneous
variety for a semisimple algebraic group decomposes in the category of Chow motives into Tate twists
of other projective homogeneous varieties we have:
\end{emptythm}

\begin{emptythm}
\label{DecompositionThm}
{\bf Theorem.}
{\it
Let~$\BM_{\ast}$ be any oriented cohomology theory (in the sense of \cite[Def. 1.1.2]{AC}),
$k$~a field, and~$X$ an isotropic projective homogeneous
$k$-variety for a semisimple algebraic group~$G$ over the field~$k$. Then there exist projective
homogeneous $k$-varieties $Y_{1},\ldots ,Y_{l}$, $l\geq 2$,
and integers $n_{i}\geq 0$, $1\leq i\leq l$, such that
$$
X\,\simeq\,\bigoplus\limits_{i=1}^{l} Y_{i}\{n_{i}\}
$$
in $\Mot_{\BM}(k)$.
}
\end{emptythm}

\goodbreak
\section{Proof of the Rost nilpotence principle for generically constant coherent theories and projective homogeneous varieties}
\label{PfRNSect}\bigbreak

\begin{emptythm}
\label{SubLemSubSect}
{\bf A technical lemma.}
Most proofs of the Rost nilpotence principle for a variety in the category of
Chow motives use a lemma proven by Rost in~\cite{Ro98}, the so called {\it Rost lemma}.
(An exception is the argument of the second author for projective quadrics.)
We prove here the analog of this lemma for a coherent theory~$\BM_{\ast}$ following
Brosnan's~\cite{Bro03} arguments for $\BM_{\ast}=\CH_{\ast}$. However there are
some subtle technical difficulites, which are addressed by the following lemma. We overcome them
using resolutions of singularities, \ie by taking advantage of our assumption that the base
field~$k$ has characteristic~$0$.

\medbreak

\noindent
{\bf Lemma.}
{\it
Let $\BM_{\ast}$ be a coherent theory over the field~$k$ and
$X,W$ be smooth schemes over~$k$ with $X$ projective. Assume~$W$ is integral with
generic point~$w$. Denote by~$\iota$ the inclusion of the gernic point $\Spec k(w)\hookrightarrow W$.
Let $\alpha\in\End_{k}(X)_{\BM}$ and $\gamma\in\BM_{i}(W\times X)$ for some~$i\in\Z$. If
$$
\alpha_{k(w)\,\ast}\big((\iota\times\id_{X})^{\ast}(\gamma)\big)\, =\; 0
$$
then there exists a closed subscheme $j:Z\hookrightarrow W$, which is $\not= W$, such
that~$\alpha_{\ast}(\gamma)$ is in the image of
$$
(j\times\id_{X})_{\ast}\, :\;\BM_{i}(Z\times X)\,\too\,\BM_{i}(W\times X)\, .
$$
}

\begin{proof}
We give a comprehensive proof indicating all properties of coherent theories which we use.

\medbreak

We show first
\begin{equation}
\label{RLpfEq1}
\alpha_{k(w)\,\ast}\big((\iota\times\id_{X})^{\ast}(\gamma)\big)\, =\,
(\iota\times\id_{X})^{\ast}\big(\alpha_{\ast}(\gamma)\,\big)\, .
\end{equation}
To prove this equation we use the following two commutative diagrams,
where we have set $\times=\times_{k}$:
\begin{equation}
\label{RLpfDiag1}
\xymatrix{
 & & X\times X &
 \\
k(w)\times X & k(w)\times X\times X \ar[l]_-{q_{13}} \ar[ru]^-{ q} \ar[d]_-{q_{12}} \ar[rr]_-{\iota\times\id_{X\times X}} & &
    W\times X\times X \ar[lu]_-{p} \ar[d]^-{p_{12}}
\\
 & k(w)\times X \ar[rr]_-{\iota\times\id_{X}} & & W\times X
}
\end{equation}
and
\begin{equation}
\label{RLpfDiag2}
\xymatrix{
k(w)\times X\times X \ar[rr]^-{\iota\times\id_{X\times X}} \ar[d]_-{q_{13}} & &
     W\times X\times X \ar[d]^-{p_{13}}
\\
k(w)\times X \ar[rr]_-{\iota\times\id_{X}} & & W\times X \rlap{\, ,}
}
\end{equation}
where~$p_{ij}$ and~$q_{ij}$ denote the projections to the respective $ij$-components.

\smallbreak

We compute now:

\smallbreak

\noindent
$\alpha_{k(w)\,\ast}\big( (\iota\times\id_{X})^{\ast}(\gamma)\big)$

\medbreak

$$
\begin{array}{c@{\;}l@{\quad}l}
          = & q_{13\,\ast}\big(\alpha_{k(w)}\cdot q_{12}^{\ast}((\iota\times\id_{X})^{\ast}(\gamma))\big) &
                \mbox{by def.} \\[3mm]
          = & q_{13\,\ast}\big( (\iota\times\id_{X\times X})^{\ast}(p^{\ast}(\alpha))\cdot
                  (\iota\times\id_{X\times X})^{\ast}(p_{12}^{\ast}(\gamma))\big)
             & \mbox{by def., Diag.\ (\ref{RLpfDiag1}), Lem.\ \ref{NatPropLem}~(ii)} \\[3mm]
          = & q_{13\,\ast}\Big( (\iota\times\id_{X\times X})^{\ast}\big(p^{\ast}(\alpha)\cdot
                p_{12}^{\ast}(\gamma)\big)\Big) & \mbox{by Lem.\ \ref{NatPropLem}~(i)} \\[3mm]
          = & (\iota\times\id_{X})^{\ast}\Big( p_{13\,\ast}\big(p^{\ast}(\alpha)\cdot
                p_{12}^{\ast}(\gamma)\big)\Big) & \mbox{by Diag.\ (\ref{RLpfDiag2}), Lem.\ \ref{NatPropLem}~(iii)} \\[3mm]
          = & (\iota\times\id_{X})^{\ast}\big(\alpha_{\ast}(\gamma)\big) & \mbox{by definition.}
\end{array}
$$
Equation~(\ref{RLpfEq1}) implies the lemma: since $\BM_{\ast}$ is coherent, $\BM_{\ast}(k(w)\times X)$ is isomorphic
to the direct limit of all $\BM_{\ast}(U\times X)$, where~$U$ runs through the open subschemes
of~$W$. Hence the assumption $\alpha_{k(w)\,\ast}\big((\iota\times\id_{X})^{\ast}(\gamma)\big)=0$
and~(\ref{RLpfEq1}) imply the existence of an open $U\subseteq W$, such that
$\alpha_{\ast}(\gamma)|_{U\times X}=0$. Let $j:Z\hookrightarrow W$ be the closed complement
of~$U$. We have $Z\not= W$, and by the localization sequence
$$
\BM_{i}(Z\times X)\,\xrightarrow{\; (j\times\id_{X})_{\ast}\;}\,\BM_{i}(W\times X)\,
   \too\,\BM_{i}(U\times X)\,\too\, 0
$$
there exists~$\delta\in\BM_{i}(Z\times X)$, such that
$\alpha_{\ast}(\gamma)=(j\times\id_{X})_{\ast}(\delta)$. We are done.
\end{proof}
\end{emptythm}

\begin{emptythm}
\label{RostLem}
{\bf The Rost lemma for coherent theories.}
{\it
Let~$\BM_{\ast}$ be a coherent theory over the field~$k$, $X,Y$ smooth and projective $k$-scheme,
and~$\alpha\in\End_{k}(X)_{\BM}$. If
$$
\alpha_{k(y)\,\ast}\big(\BM_{\ast}(X_{k(y)})\big)\, =\, 0
$$
for all~$y\in Y$ then
$$
(\alpha^{\circ (1+\dim Y)})_{\ast}\big(\BM_{\ast}(Y\times_{k}X)\big)\, =\, 0.
$$
}

\begin{proof}
Let~$n=\dim Y$.

\smallbreak

We define a filtration on $\BM_{\ast}(Y\times_{k}X)$ as follows. Let~$F_{l}\subseteq\BM_{\ast}(Y\times_{k}X)$
for all $0\leq l\leq n$ be the subgroup generated by all images
$$
(j_V\times\id_{X})_{\ast}\, :\;\BM_{\ast}(V\times_{k} X)\,\too\,
\BM_{\ast}(Y\times_{k}X)\, ,
$$
where $j_V:V\too Y$ is a closed $k$-subscheme of~$Y$ of dimension~$\leq l$, equivalently
we can take only the images with~$V$ a closed integral $k$-subscheme.

\smallbreak

We have
$$
\BM_{\ast}(Y\times_{k}X))\, =\, F_{n}\,\supseteq\, F_{n-1}\,\ldots\,\supseteq\, F_{1}\,
\supseteq\, F_{0}\,\supseteq\, F_{-1}\, :=\;\{ 0\}\, ,
$$
and so it is enough to show that $\alpha_{\ast}(F_{l})\subseteq F_{l-1}$
for all $-1\leq l\leq n:=\dim Y$. We verify this by induction on~$l$.
The induction beginning $l=-1$ is clear, so let $l\geq 0$.

\medbreak

Let for this~$a\in F_{l}\subseteq\BM_{\ast}(Y\times_{k}X)$ be equal to $(j_{V}\times\id_{X})_{\ast}(b)$, where
$$
(j_{V}\times\id_{X})_{\ast}\, :\;\BM_{\ast}(V\times X)\,\too\,\BM_{\ast}(Y\times_{k}X)\, ,
$$
and $j_{V}:V\hookrightarrow Y$ is a closed embedding of a subscheme of dimension $\leq l$.
We have to show $\alpha_{\ast}(a)\in F_{l-1}$, and can assume for this that~$V$ is a
closed integral subscheme.

\smallbreak

Let $\pi:W\too V$ be a projective birational $k$-morphism with~$W$ smooth projective and integral,
which exists since~$k$ has characteristic~$0$ by Hironaka~\cite{Hi}. By Lemma~\ref{birOntoLem} there
exists then $c\in\BM_{\ast}(W\times_{k}X)$, a closed subscheme $j:V'\hookrightarrow V$,
and $b'\in\BM_{\ast}(V'\times_{k}X)$, such that
$$
b\, =\, (\pi\times\id_{X})_{\ast}(c)+(j\times\id_{X})_{\ast}(b')\, .
$$
Hence
$$
a\, =\,( \varepsilon\times\id_{X})_{\ast}(c)+
( j_{V'}\times\id_{X})_{\ast}(b'),
$$
where $\varepsilon=j_{V}\circ\pi$ and $j_{V'}=j_V\circ j$.
Note now that $\dim V'<\dim V\leq l$ since~$V$ is irreducible, and so we have
that $\big( j_{V'}\times\id_{X}\big)_{\ast}(b')\in F_{l-1}$,
which by induction implies
$$
\alpha_{\ast}\big(\,( j_{V'}\times\id_{X})_{\ast}(b')\,\big)
\;\subseteq\, F_{l-2}\subseteq F_{l-1}\, .
$$
Hence it is enough to show that
$\alpha_{\ast}\big(( \varepsilon\times\id_{X})_{\ast}(c)\big)$ is
in~$F_{l-1}$.

\medbreak

It follows from the standard compatibilities of pull backs and push forwards that the diagram
$$
\xymatrix{
{\BM_{\ast}}(W\times_{k}X) \ar[rr]^-{(\varepsilon\times\id_{X})_{\ast}} \ar[d]_-{\alpha_{\ast}} & &
    {\BM_{\ast}}(Y\times_{k}X) \ar[d]^-{\alpha_{\ast}}
\\
{\BM_{\ast}}(W\times_{k}X) \ar[rr]_-{(\varepsilon\times\id_{X})_{\ast}} & & {\BM_{\ast}}(Y\times_{k}X)
}
$$
is commutative. Hence
$\alpha_{\ast}\big(( \varepsilon\times\id_{X})_{\ast}(c)\big)=( \varepsilon\times\id_{X})_{\ast}(\alpha_{\ast}(c))$.

Let $v,w$ be the generic points of~$V$ and~$W$, respectively. We have $k(v)=k(w)$
and so by our assumption
$$
\alpha_{k(w)\,\ast}\big(\, (\iota\times\id_{X})^{\ast}(c)\big)\, =\, 0\, ,
$$
where $\iota:\Spec k(w)\hookrightarrow W$ is the embedding of the generic point.
Therefore by the Lemma in~\ref{SubLemSubSect} there exists a closed
embedding $q:Z\hookrightarrow W$ with $Z\not= W$ and so $\dim Z<\dim W\leq l$,
such that $\alpha_{\ast}(c)$ is in the image of
$$
(q\times\id_{X})_{\ast}\, :\;\BM_{\ast}(Z\times_{k}X)\,\too\,\BM_{\ast}(W\times_{k}X)\, .
$$
Then $\alpha_{\ast}\big(( \varepsilon\times\id_{X})_{\ast}(c)\big)$ is in the image of $(p\times\id_{X})_{\ast}$, where
$p=\varepsilon\circ q:Z\rightarrow Y$. And hence, in the image of $(j_{\overline{Z}}\times\id_{X})_{\ast}$, where
$\overline{Z}$ is the reduced image of $p$ and $j_{\overline{Z}}:\overline{Z}\rightarrow Y$ is the respective closed embedding.
Since $\dim \overline{Z}<l$, this element belongs to $F_{l-1}$.
\end{proof}
\end{emptythm}

\begin{emptythm}
\label{RostL-wrong-non-coh}
{\bf Remark.}
For non-coherent theories, the Rost Lemma is in general wrong. Here is an example.
Let $\BM_{\ast}$ be $\operatorname{CH}^*_{alg}$, and $X=Y$ be an elliptic curve $C$.
Let $\alpha=[\Delta]-[C\times p]-[p\times C]$, where $\Delta$ is the diagonal and $p$ is some
$k$-rational point on $C$. Then, for any $L/k$, $\operatorname{CH}^*_{alg}(X_L)={\Bbb{Z}}\cdot 1\oplus {\Bbb{Z}}\cdot[p]$,
and so ${\alpha_L}_*(\BM_{\ast}(X_L))=0$. At the same time,
$\alpha$ is a {\it projector} $\alpha^{\circ 2}=\alpha$. Hence, for any natural $N$,
$\alpha^{\circ N}_*(\BM_{\ast}(Y\times X))\neq 0$, since $\alpha\neq 0$.

\medbreak

Using the decomposition from Theorem~\ref{DecompositionThm} the Rost lemma for~$\BM_{\ast}$
implies the Rost nilpotence principle for projective homogeneous varieties and generically constant coherent theories.
The proof is word by word the same as~\cite[Proof of Thm.\ 8.2]{ChGiMe05} (only replacing
$\CH_{\ast}$ by~$\BM_{\ast}$). However to indicate where we need that the product in a coherent theory
is compatible with base change we give the details.
\end{emptythm}

\begin{emptythm}
\label{RNThm}
{\bf Corollary.}
{\it
Let~$X$ be a projective homogeneous variety for a semisimple algebraic group
over the field~$k$. Then~$X$ satisfies the Rost nilpotence principle for every
generically constant coherent theory~$\BM_{\ast}$.
}

\begin{proof}
If~$X$ is split then it is a cellular variety and so by~\cite[Cor.\ 2.9]{ViYa07}
a direct sum of twists of Tate motives. Since $\BM_{\ast}$ is {\it generically constant},
$\BM_{\ast}(k)\too\BM_{\ast}(L)$ is an isomorphism for all field extensions~$L\supseteq k$,
and so Rost nilpotence holds trivially in this case.

\smallbreak

If~$X$ is not split then it still has a decomposition as in Theorem~\ref{DecompositionThm}:
$$
X\,\simeq\,\bigoplus\limits_{i=1}^{l}Y_{i}\{n_{i}\}
$$
for projective homogeneous varieties~$Y_{i}$ and integers~$n_{i}\geq 0$
(it certainly has such a decomposition with~$l=1$ and $Y_1=X$, $n_1=0$).
We prove now by down going induction on the number of components~$l$ that
there exists an integer~$m\geq 0$ which only depends on the dimension of the projective
homogenous variety and on the number~$l$, such that $\alpha^{\circ\, m}=0$ if~$\alpha$
is in the kernel of
$$
\End_{k}(X)_{\BM}\,\too\,\End_{L}(X_{L})_{\BM}
$$
for some field extension $L\supseteq k$.
As the number
of possible components is bounded by the number of Tate motives appearing in
a decomposition over the algebraic closure of~$k$, we have the base of induction.

\smallbreak

Let~$\alpha\in\End_{k}(X)_{\BM}$, such that $\alpha_{L}=0$ for some field
extension~$L$ of~$k$. If~$y$ is a point of~$Y_{i}$ then~$Y_{i}$ is isotropic over~$k(y)$
and so by Theorem~\ref{DecompositionThm} the motive $Y_{i\, k(y)}\{n_{i}\}$
is a direct sum of at least two Tate twisted motives of projective homogenous
varieties over~$k(y)$. Hence over~$k(y)$ the motive of~$X$ decomposes in
more than~$l$ summands and so by induction we have $\alpha_{k(y)}^{\circ\, m}=0$
for an integer~$m\geq 1$ which depends only on the number of components of
the decomposition of~$X$ over~$k(y)$ and the dimension of~$X$. As the number
of possible components is bounded we can choose~$m$, such
that $\alpha_{k(y)}^{\circ\, m}=0$ for all $y\in Y_{i}$ and all $1\leq i\leq l$.

\smallbreak

By Lemma~\ref{NatPropLem}~(iv) we have
$$
0\, =\,\alpha_{k(y)}^{\circ\, m}\, =\, (\alpha^{\circ\, m})_{k(y)}\, ,
$$
and so $(\alpha^{\circ\, m})_{k(y)\,\ast}\big(\BM_{\ast}(X_{k(y)})\big)=0$ for
all $y\in Y_{i}$ and all $1\leq i\leq l$. By the Rost lemma this implies
$$
(\alpha^{\circ\, m\cdot (1+\dim Y_{i})})_{\ast}\big(\,\Hom_{k}(Y_{i}\{n_{i}\},X)_{\BM}\,\big)\; =\, 0
$$
for all $1\leq i\leq l$. Let~$t$ be the maximum of the dimensions of the~$Y_{i}$.
Then since the motive of $X$ is equal the direct sum of the motives $Y_{i}\{n_{i}\}$'s in the category of $\BM$-motives
we get $\alpha^{\circ\, m\cdot (1+t)}=0$. We are done.
\end{proof}
\end{emptythm}

\medbreak

From \cite[Cor.\ 4.4.3]{AC} we know that {\it free} theories are generically constant and from Corollary in \ref{free-theories}
these are also coherent. Hence, we obtain:

\begin{emptythm}
\label{Rost-free}
{\bf Corollary.}
{\it
Let~$X$ be a projective homogeneous variety for a semisimple algebraic group
over the field~$k$. Then~$X$ satisfies the Rost nilpotence principle for every
free theory~$\BM_{\ast}$.
}

\medbreak

Note that not all {\it generically constant coherent} theories are {\it free} (= of {\it rational type} \cite{SU}). Here is a large supply of examples.
\end{emptythm}

\medbreak
\begin{emptythm}
\label{coherent+constant-not-free}
{\bf Example.}
Let $A^{\ast}$ be a generically constant coherent theory over $k$,
and $\{(Q_{\lambda},a_{\lambda})\}_{\lambda\in\Lambda}$ a
collection of smooth projective varieties $Q_{\lambda}$ defined over $k$ (!) with some classes
$a_{\lambda}\in A^{\ast}(Q_{\lambda})$. In other words, we have a collection of $A^*$-correspondences
$\rho_{\lambda}:Q_{\lambda}\rightsquigarrow\Spec(k)$. Suppose that the map $\rho_*: A_*(Q_{\lambda,F})\rightarrow A_*(F)$
sending $v$ to $\pi_*(v\cdot a_{\lambda})$ is zero, for any field extension $F/k$.
For any $F/k$, and any smooth variety $X/F$ we have natural projections:
$$
\xymatrix{
Q_{\lambda,F} & Q_{\lambda,F}\times_F X \ar[l]_(0.6){\pi_1} \ar[r]^(0.65){\pi_2} & X
}
$$

Now, in $A^{\ast}(X)$ we can mod-out the images of the maps from $\rho_*:A^{\ast}(Q_{\lambda,F}\times_F X)\rightarrow A^{\ast}(X)$ sending
$u$ to $(\pi_2)_{\ast}(u\cdot\pi_1^{\ast}(a_{\lambda}))$.
It is not difficult to check that the resulting theory $\widetilde{A}^{\ast}$ is coherent and generically constant.
As long as not all $a_{\lambda}$'s are zero, it will not be free, since it has the same coefficient ring as $A^{\ast}$,
but the natural projection $A^{\ast}\rightarrow\widetilde{A}^{\ast}$ has a non-trivial kernel, and so, there is a non-trivial
kernel of the projection $\widetilde{A}^{\ast}(k)\otimes_{\Laz}\Omega^{\ast}\twoheadrightarrow\widetilde{A}^{\ast}$.

\smallbreak

Here is the simplest such example. Take $A^{\ast}=\operatorname{CH}^{\ast}$ and a single elliptic curve $C$ over $k$ with $a=[p_1]-[p_0]$,
where $p_0,p_1$ are two distinct $k$-rational points on $C$. Then the resulting theory $\widetilde{A}^{\ast}$ will be
a constant coherent theory in-between $\operatorname{CH}^{\ast}$ and $\operatorname{CH}^{\ast}_{alg}$.
\end{emptythm}

\bigbreak

\begin{emptythm}
\label{failure-Rost}
{\bf Failure of the Rost nilpotence principle for non-coherent theories.}\
For non-coherent theories it does not make much sense to speak about Rost nilpotence principle.
If the theory is non-constant, counterexamples are obvious: just take $A^*_k=\operatorname{CH^*}$, and $A^*_F=\operatorname{CH}^*/2$,
for any non-trivial extension $F/k$. But even for constant ones we have:

\medbreak

\noindent
{\bf Example.} Let $C$ be some elliptic curve over $k$, and $p$ - a $k$-rational point on it.
Take $B^*_k=\operatorname{CH}^*/2$. Consider the $B^*$-correspondence $\rho: C\times C\rightsquigarrow\Spec(k)$
given by the class $\alpha=[\Delta]-[C\times p]-[p\times C]$. Notice, that it has the property, that, for all $F/k$,
the map $\rho_*:\operatorname{CH}^*/2((C\times C)_F)\rightarrow\operatorname{CH}^*/2(F)$ is zero.

\smallbreak

For a non-trivial extension $F/k$, and smooth $X/F$ let us take $B^*_F(X)$ to be the quotient of $\operatorname{CH}^*/2(X)$ modulo the
image of $\rho_*:\operatorname{CH}^*/2(C\times C\times X)\rightarrow\operatorname{CH}^*/2(X)$.
Then $B^*$ will be an extended constant oriented cohomology theory, and we have a $B^*$-correspondence
$\alpha:C\rightsquigarrow C$ which vanishes over all non-trivial field extensions by construction.
At the same time, $\alpha$ itself is a projector. This projector is clearly non-trivial, since the motive of an elliptic curve
(even with $\Z/2$-coefficients) is not a direct sum of two Tate-motives.
Hence, the Rost nilpotence principle fails for the theory $B^*$.
\end{emptythm}

\bigbreak

\bibliographystyle{amsalpha}

\end{document}